\documentclass{article}
\usepackage{amssymb}

\usepackage{graphicx}
\usepackage{amsmath}
\usepackage{a4wide}

\numberwithin{equation}{section}

\newcommand{\Ad}{\operatorname{Ad}}

\input{tcilatex}

\begin{document}

\begin{center}
\bigskip {\Large NAHM'S EQUATIONS, CONFIGURATION SPACES AND FLAG MANIFOLDS}

by\bigskip

MICHAEL ATIYAH and ROGER BIELAWSKI$^1$\bigskip
\end{center}

\addtocounter{footnote}{1}\footnotetext{EPSRC Advanced Fellow}

\section{Introduction}

This paper has an unusual origin, evolution and potential application. \ As
explained in \cite{At1, At2} it arose from a problem posed by Berry and
Robbins \cite{BR} in their investigation of the spin-statistics theorem. \
They asked the following simple question: is there, for each integer $n\geq
2,$ a continuous map
\begin{equation}
f_{n}:C_{n}({\Bbb R}^{3})\rightarrow U(n)/T^{n}  \label{1.1}
\end{equation}
compatible with the action of the symmetric group $\Sigma _{n}?$ \ Here $%
C_{n}({\Bbb R}^{3})$ is the configuration space of $n$ ordered distinct points of $%
{\Bbb R}^{3},$ and $U(n)/T^{n}$ is the well known flag manifold. \ The symmetric
group acts freely on both spaces, by permuting points in the first space and
components of the flag in the second. \ For $n=2,$%
\begin{eqnarray*}
C_{2}({\Bbb R}^{3}) &=&{\Bbb R}^{3}\times ({\Bbb R}^{3}-0) \\ U(2)/T^{2} &=&P_{1}({\Bbb C})=S^{2}
\end{eqnarray*}
and there is an obvious solution to (1.1). \ Note that this obvious solution
is also compatible with the natural action of $SO(3)$ on both sides.

In \cite{At1} a positive answer was given to the Berry-Robbins question using
an elementary construction of $f_{n}.$ \ A more elegant construction was also
proposed but this was dependent on the conjectured non-vanishing of a certain
determinant. \ This question was pursued further in \cite{At2} and the
conjecture has now been verified numerically for $n\leq 20$ \cite{AtSut}.

The maps $f_{n}$ of \cite{At1} are all compatible with the action of $SO(3),$
where we choose $SO(3)$ to act on $U(n)/T^{n}$ via its irreducible
representation on ${\Bbb C}^{n}.$ \ This suggested a natural generalization of the
Berry-Robbins question to other compact Lie groups $G$ instead of $U(n).$ \
Let $T$ be a maximal torus of $G,$ then the Weyl group
\begin{equation*}
W=N(T)/T
\end{equation*}
acts freely on the flag manifold $G/T.$ \ Let $\frak{h}$ be the Lie algebra
of $T,$ then $W$ acts also on $\frak{h}$ and on
\begin{equation}
\frak{h}^{3}=\frak{h}\otimes {\Bbb R}^{3}.  \label{1.2}
\end{equation}
The singular set $\Delta $ of this action on $\frak{h}^{3}$ is the union of
the codimension $3$ subspaces which are the kernels of root homomorphisms
\begin{equation*}
\alpha\otimes 1 :\frak{h}^{3}\rightarrow {\Bbb R}^{3.}
\end{equation*}
Then $W$ acts freely on\textit{\ }$\frak{h}^{3}-\Delta $ which is
the space of \textbf{regular }triples in $\frak{h}$ (i.e. with only $\frak{h}
$ as their common centralizer).

For $G=U(n)$ we recognise that $G/T$ is the usual flag manifold and that
\begin{equation*}
\frak{h}^{3}-\Delta =C_{n}({\Bbb R}^{3}).
\end{equation*}
The obvious generalization of (1.1) is therefore to ask for a continuous map
\begin{equation}
f_{G}:\frak{h}^{3}-\Delta \rightarrow G/T  \label{1.3}
\end{equation}
which is compatible with the action of the Weyl group. \ Again we can hope
to find $f_{G}$ which is also compatible with the action of $SO(3)$, where $%
SO(3)$ acts on $h^{3}-\Delta ,$ via the decomposition (1.2), and acts on $%
G/T $ through some preferred homomorphism
\begin{equation}
\rho :SU(2)\rightarrow G.  \label{1.4}
\end{equation}
There is a natural candidate for each compact Lie group $G,$ generalizing the
irreducible $n$-dimensional representation of $SU(2)$ for $U(n).$ \ This
is given by the so-called \textbf{regular} (or \textbf{principal) }%
homomorphism $\rho .$ \ This may be characterized by the fact that, after
complexification, $\rho $ takes the unipotent element of $SL(2, {\Bbb C})$ into a
\textbf{regular } unipotent element of $G^{\Bbb C}(i.e.$ one which lies in a
unique Borel subgroup). \ The regular homomorphism is unique\ up to
conjugacy. \ Its action on $G/T$ also factors through $SO(3).$

It turns out that such a map $f_{G},$ with all the desired properties, can
actually be extracted from previous work on \textbf{Nahm's equations }%
  in
\cite{BielCMP}. \ The original purpose of this paper was to show how this
comes about.

After the original solution of the Berry-Robbins problem in \cite{At1}, various
cohomological consequences were drawn in \cite{At3}, and similar results were
expected for other Lie groups. \ It was then suggested by Gus Lehrer that
these ideas might be related to the Springer representation of the Weyl group
and the extensive work done in this direction by Kazhdan and Lusztig (see for
example \cite{KL} or \cite{KL1}). \ This has led us to extend our
investigations, using Nahm's equation, to include arbitrary homomorphisms $%
\rho $ of $SU(2)$ into $G.$ \ This leads to an interesting geometrical
picture, generalizing the map (1.3). \ It is our hope that this will shed
light on the work of Kazhdan and Lusztig and explain the geometry behind the
Hecke algebras.

The paper is organised as follows. \ In \S 2 we review the key aspects of
Nahm's equations and the various moduli spaces. \ Then in \S 3 we spell out
the main construction which in particular gives the map (1.3). \ In \S 4 we
break the $SU(2)$-symmetry down to a circle subgroup and relate the geometry
to that of the complex Lie group. \ In \S 5 we explain the relation of our
construction to the Kazhdan-Lusztig work.

In order to keep the geometrical picture clear Sections 2-5 are presented in
non-technical terms. \ The precise analytical details are then set out in section 6.

\section{Nahm's Equations and Lie Groups \label{2}}

Since Nahm's equations will be our main technical tool it may be helpful to
provide here a little background on how these equations first arose and what
role in particular they play in Lie theory.

For any Lie algebra $\frak{g}$ Nahm's equations are the system of $3$ 
${\frak g}$-valued
ordinary differential equations
\begin{equation}
\frac{dT_i}{dt}+[{T}_{j},T_{k}]=0  \label{2.1}
\end{equation}
when $(i,j,k)$ is a cyclic permutation of $(1,2,3)$ and the $T_{i}$ are
functions of the real variable $t.$

While (2.1) makes sense for any Lie algebra these equations have a
particularly simple interpretation when $\frak{g}$ is the Lie algebra \ of a
\textbf{compact }Lie group $G$ (the case of interest to us). \ In this case
we have a $G$-invariant metric $\langle ,\rangle $ on $\frak{g},$ enabling us
to
identify $\frak{g}$ with its dual $\frak{g}^{\ast }.$ \ This leads to the well-known $G$%
-invariant skew $3$-form $\phi $ on $\frak{g}$ given by
\begin{equation}
\phi (T_{1},T_{2},T_{3})=\langle T_{1},[T_{2},T_{3}]\rangle  \label{2.2}
\end{equation}
which also defines the bi-invariant (harmonic) exterior differential $3$%
-form on $G$ (unique up to a scalar for simple $G).$ \ It is then easy to
check that
\begin{equation}
\begin{tabular}{l}
$\text{\textbf{Nahm's equations are the gradient-flow }}$ \\
$\text{\textbf{equations for} }\phi\quad \text{\textbf{as a function on
}}\frak{g}\oplus
\frak{g}\oplus \frak{g}.$%
\end{tabular}
\label{2.3}
\end{equation}

Regarding $\phi $ as a function on $\frak{g}\otimes {\Bbb R}^{3}$ or on $\func{Hom}%
({\Bbb R}^{3},\frak{g})$ it is \textbf{invariant under the }$SO(3)$-\textbf{action on }$%
{\Bbb R}^{3}.$

To see this we observe that if
\begin{equation*}
T:{\Bbb R}^{3}\rightarrow \frak{g}
\end{equation*}
is a linear map then
\begin{equation*}
\Lambda ^{3}T:\Lambda ^{3}{\Bbb R}^{3}\rightarrow \Lambda ^{3}\frak{g}
\end{equation*}
sends the $SO(3)$-invariant oriented volume element of ${\Bbb R}^{3}$ to $T_1\wedge T_2\wedge T_3 $.

In Lie theory it is standard to consider
\begin{equation*}\frak{g}\oplus \frak{g}\cong \frak{g}\otimes {\Bbb R}^{2}\cong \frak{g}^{\Bbb C}
\end{equation*}
the complexified Lie algebra. \ What, one may ask, is the significance of
replacing $\frak{g}\otimes {\Bbb R}^{2}$ by $\frak{g}\otimes {\Bbb R}^{3}$ as in Nahm's
equations?
\
The answer is that we should identify ${\Bbb R}^{3}$ here with the \textbf{%
imaginary quaternions}
\begin{equation*}
{\Bbb R}^{3}\cong \func{Im}({\Bbb H}).
\end{equation*}
To see why this is the case we should actually introduce a fourth
Lie-algebra-valued function $T_{0}(t)$ and consider the expression
\begin{equation}
A=T_{0}dt+T_{1}dx_{1}+T_{2}dx_{2}+T_{3}dx_{3}  \label{2.4}
\end{equation}
as defining a $G$-connection over
\begin{equation*}
H={\Bbb R}^{4}={\Bbb R}\oplus \func{Im}{\Bbb H}.
\end{equation*}
Since the matrices in (2.4) depend only on $t,$ and not on $x,$ the natural
gauge group to consider is simply the $G$-valued functions of $t.$ \ Using
the gauge freedom we can reduce $T_0$ to zero getting back to just $3$
matrices $T_{1},T_{2},T_{3}.$ \ More invariantly we should start with a $G$%
-bundle over ${\Bbb R}^{4}$ which has an action of the translations of ${\Bbb R}^{3}.$ \
Then (2.4) describes a connection for this bundle which is ${\Bbb R}^{3}$%
-invariant and is written in an ${\Bbb R}^{3}$-invariant gauge. \ The matrices $%
T_{i} $ then represent the difference between the Lie derivative and the
covariant derivative in the $i^{\text{th}}$ direction, and are usually
referred to as Higgs fields: they are infinitesimal automorphisms of the
bundle.

Now in $4$ dimensions we have the famous anti-self-duality (ASD) equations
\begin{equation*}
\ast F=-F
\end{equation*}
where $F$ is the curvature of a connection $A.$ \ It was Donaldson \cite{Don}
who first observed that, for the connection (2.4), and after gauging away $%
T_{0},$ \textbf{the }ASD \textbf{equations are identical with Nahm's
equations.}

Now it is an important point that the ASD equations over ${\Bbb R}^{4}$ are,
formally, the \textbf{hyperk\"{a}hler moment map} for the action of the gauge
group. \ This leads (formally) to a hyperk\"{a}hler metric on moduli spaces
of solutions. \ This observation is well-known to physicists as a consequence
of super-symmetry and the concept of the hyperk\"{a}hler quotient
construction, developed in \cite{HKLR}, was inspired by this.

Hyperk\"{a}hler manifolds are Riemannian manifolds of dimension $4n$ with
holonomy in $Sp(n),$ so that their tangent spaces are quaternionic. \ They
have a $2$-parameter family of complex structures (depending on an embedding
${\Bbb C}\rightarrow {\Bbb H},$ or on an imaginary quaternion $I$ with $I^{2}=-1).$ \ They
are the quaternionic counterparts of complex K\"{a}hler manifolds and they
have twistor spaces in the sense of Roger Penrose.

All these general remarks apply, not only to the full four-dimensional ASD
equations (where the matrices $T_{0},...,T_{4}$ in (2.4) depend on all $4$
variables), but also to the (partially) translation invariant ones such as
Nahm's equations (where the dependence is only on one variable). \ These
moduli spaces of solutions to Nahm's equations should formally have
hyperk\"{a}hler metrics.

Of course appropriate boundary conditions need to be imposed and analytical
details need to be checked. \ Originally Nahm introduced his equations in
relation to non-abelian magnetic monopoles (which satisfy the Bogomolny
equations, the ${\Bbb R}$-invariant version of the ASD equations) and the
corresponding hyperk\"{a}hler metrics were studied in detail in \cite{AtHi}.

It was Kronheimer \cite{Kron, Kron2} who first applied Nahm's equations to the
study of Lie groups themselves, by altering the boundary conditions. \ For
$SU(2)$-monopoles of charge $n$ Nahm considered his equation for $G=U(n)$ on an interval and took
as boundary condition that the $T_{i}$ had \textbf{simple regular poles} at
each end. \ If the $T_{i}$ have simple poles with residues $\sigma _{i}$ then
(2.1) shows that
\begin{equation}
\sigma _{i}=[\sigma _{j},\sigma _{k}]  \label{2.5}
\end{equation}
are the commutation relations (up to a factor $2)$ of the quaternions $%
i ,j,k$%
\begin{equation*}
\sigma _{1}=i/2,\;\sigma _{2}=j/2,\;\sigma _{3}=k/2
\end{equation*}
and thus are the standard generators of the Lie algebra of $SU(2).$ \ A pole
is called regular if the $n$-dimensional representation of $\frak{su}(2)$
given by the matrices (2.5) is irreducible.

In \cite{Kron2} Kronheimer considered poles of any type, characterized by an arbitrary
homomorphism
\begin{equation*}
\rho :\frak{su}(2)\rightarrow \frak{g} \label{2.6}
\end{equation*}
given by the residues as in \eqref{2.5}. \ By taking $\rho =0$ at one end (so that
the solution has no pole there), and a general $\rho $ at the other,
Kronheimer obtained as his moduli space a new hyperk\"{a}hler manifold which
(for almost all of its complex structures) could be identified with the
nilpotent orbit in $\frak{g}^{\Bbb C}$ corresponding to $\rho $ (i.e. the one
containing $\rho (x)$ where $x$ is a nilpotent element of $\frak{sl}(2,C)$
and $\rho $ is understood as the complexification of (2.6)).

In \cite{Kron} Kronheimer considered Nahm's equation on the half-line $t\geq 0$
and imposed finiteness at $0,$ and finite limits at $\infty $%
\begin{equation*}
T_{i}(t)\rightarrow T_{i}(\infty )=\tau _{i},
\end{equation*}
where $T_{1}(\infty ),\;T_{2}(\infty ),\;T_{3}(\infty )$ are a regular
(commuting) triple. \ For these boundary conditions (with the $G$-conjugacy
class of the regular triple $\tau $ fixed) he found the moduli space to be a
hyperk\"{a}hler manifold which
(for almost all of its complex structures) was a regular semi-simple orbit in ${\frak g}^{\Bbb C}$. 

These results of Kronheimer have since \cite{Biq, Kov} been extended to provide
hyperk\"{a}hler metrics for all complex co-adjoint orbits. \ This story is
the quaternionic generalization of the complex K\"{a}hler metrics on
co-adjoint orbits of $G.$

The moral of all this is the following. \ A compact (real) Lie group $G$ has
a complexification $G^{\Bbb C}$ with compact complex homogeneous spaces (e.g. $%
G^{\Bbb C}/B)$ which have K\"{a}hler metrics. \ The Lie algebra $\frak{g}$ has a
(vector space) quaternionisation $\frak{g}\otimes {\Bbb H},$ but there is no corresponding
``quaternionic group''. \ \ However the analogous ``homogeneous spaces'' do
exist as hyperk\"{a}hler manifolds. \ For many purposes $G^{\Bbb C}$ can be
studied through for example the flag manifold $G^{\Bbb C}/B,$ so we can view the
hyperk\"{a}hler structures on the complex co-adjoint orbits of $G^{\Bbb C}$ as
substitutes for the non-existing quaternion group.

In this spirit the different homomorphisms $\rho :SU(2)\rightarrow G$ are the
quaternionic analogous of $1$-parameter subgroups $U(1)\rightarrow G.$

From this point of view Nahm's equation is the key to unlocking the
``quaternionic nature of Lie groups''. \ An area where this has proved its
worth is in the clarification of the work of Brieskorn on Kleinian
singularities (due to Kronheimer \cite{KronALE}) and its systematic extension to
the Brieskorn-Grothendieck resolution of singularities of the nilpotent
variety \cite{Slo}.

In \cite{BielJLMS, BielCMP} other variants of the boundary conditions for Nahm's equations
were studied. \ Here we shall be concerned with the equations on the
half-line where, following Kronheimer, we take limiting regular triples at $%
\infty $ but as $t\rightarrow 0$ we impose a simple pole of type $\rho .$ \
The case when $\rho $ is the regular $SU(2)$ will give the construction of
the map (1.3), while the other $\rho $ will yield the more general picture to
be discussed later.

\section{The Main Construction}

Let $\rho :su(2)\rightarrow \frak{g}$ be a homomorphism, and consider
solutions of Nahm's equations (2.1) on the half-line $0<t<\infty ,$ with the
boundary conditions:
\begin{equation}
\begin{tabular}{ll}
(a) & there is a pole of type $\rho $ as $t\rightarrow 0$ \\
(b) & the $T_{i}$ tend to a regular commuting triple in $\frak{g}^{3}$ as $%
t\rightarrow \infty .$%
\end{tabular}
\label{3.1}
\end{equation}
We denote the space of such solutions by $N^{\prime }(\rho ).$ \ By taking
the value at $\infty $ we get a map
\begin{equation}
N^{\prime }(\rho )\rightarrow \frak{g}^{3}.  \label{3.2}
\end{equation}

Now fix a maximal torus $T$ of $G$ and let $\frak{h}$ be its Lie algebra. \
$G$
acts on $\frak{g}^{3}$ and on the regular commuting triples. \ Each orbit is of the form $%
G{\tau }$ where $\tau $ is a regular triple of $\frak{h}$ and every orbit
$G{\tau }$ meets $\frak{h}^{3}$ in an orbit of the Weyl group $W.$ \ We can
therefore define a finite covering $N(\rho )$ of $N^{\prime }(\rho )$ by the
commutative diagram
\begin{equation}
\begin{array}{ccccc}
N(\rho ) &  & \rightarrow &  & N^{\prime }(\rho ) \\ &  &  &  &  \\
\downarrow &  &  &  & \downarrow \\ &  &  &  &  \\ \frak{h}^{3}-\Delta &  &
\rightarrow &  & (\frak{h}^{3}-\Delta )/W
\end{array}
\label{3.3}
\end{equation}
where the vertical arrows assign to a solution of Nahm's equation its orbit
type at $\infty ,$ arising from (3.2).

Fixing $\tau $ identifies $G{\tau }$ with $G/T$ and hence, by taking the values at $\infty$,  we get a natural
map
\begin{equation}
\phi (\rho ):N(\rho )\rightarrow G/T.  \label{3.4}
\end{equation}

$N(\rho )$ is a fibration over $\frak{h}^{3}-\Delta $ with fibre (at $\tau )$ $%
N(\rho ,\tau )$ and the manifolds $N(\rho ,\tau )$ are all hyperk\"{a}hler. \
In fact if we denote by $M(\rho )$ the $T$-bundle over $N(\rho )$ induced by
$\phi $ then $M(\rho )$ is a hyperk\"{a}hler manifold and the map
\begin{equation*}
\mu :M(\rho )\rightarrow \frak{h}^{3}-\Delta
\end{equation*}
is a hyperk\"{a}hler moment map of the $T$-action. \ The manifolds $N(\rho
,\tau )$ are just the \textbf{hyperk\"{a}hler quotients. \ }$M(\rho )$ itself
is also a suitable moduli space of solutions of Nahm's equations.

All these statements are best understood in terms of the gauged version of
Nahm's equations involving the fourth matrix $T_{0}.$ \ This, together with
the more precise description of the analytical details will be explained in
\S 6 .

The action of $W$ on $N(\rho ),$ implied by (3.3), is induced by an action of
the normalizer $N(T)$ on $M(\rho).$ \ Moreover the group $SU(2)$ acts
throughout, commuting with $N(T),$ and the map $\phi (\rho )$ of (3.4)
is compatible with the $SU(2)$ action on $G/T$ induced by $\rho .$

In fact all these constructions are compatible with yet another group. \ This
is the group $Z(\rho ),$ the centralizer of $\rho (SU(2))$ in $G.$ \
Conjugation by an element of $Z(\rho )$ preserves the boundary conditions
(3.1) and so induces an action on $N(\rho ).$ \ The natural action of $%
Z(\rho )$ on $G/T$ also commutes (by definition) with the action of $SU(2).$
\ Thus $Z(\rho )$ lifts also to an action on $M(\rho ).$

To sum up we have a diagram of maps
\begin{equation}
\begin{array}{ccccc}
M(\rho ) &  & \rightarrow &  & G \\ &  &  &  &  \\ \downarrow &  &  &  &
\downarrow \\ N(\rho ) &  & \overset{\phi }{\rightarrow } &  & G/T \\ &  &  &
&
\end{array}
\label{3.5}
\end{equation}
and a compatible action of the group
\begin{equation*}
N(T)\times SU(2)\times Z(\rho )
\end{equation*}
descending to an action of
\begin{equation*}
W\times SU(2)\times Z(\rho )
\end{equation*}
for the bottom map $\phi .$ \ For the fibre map
\begin{equation}
N(\rho )\rightarrow \frak{h}^{3}-\Delta  \label{3.6}
\end{equation}
$W\times SU(2)$ acts naturally on the base, while $Z(\rho )$ acts trivially
on the base but acts on the fibres $N(\rho ,\tau ).$

The torus $T$ and the group $Z(\rho )$ both preserve the hyperk\"{a}hler
structure of $M(\rho ),$ but $SU(2)$ rotates the complex structures.

There are three noteworthy special cases of $\rho .$ \ These are

\begin{enumerate}
\item[(a)]  $\rho =0.$ \ Then $Z(\rho )=G$ and, as will be discussed in the
next section, $N(\rho ,\tau )$ is the complexification $G^{\Bbb C}/T^{\Bbb C}$ of $G/T.
$ \ The map
\begin{equation*}
\phi :G^{\Bbb C}/T^{\Bbb C}\rightarrow G/T
\end{equation*}

commutes with $G.$ \ Observing that $G/T$ sits inside $G^{\Bbb C}/T^{\Bbb C}$ with a
contractible $T$-invariant slice it follows that $\phi $ must be a
deformation retraction compatible with this $G$-action.

\item[(b)]  $\rho $ the regular $SU(2).$ \ Then $Z(\rho )$ is finite and, as will be
shown in the next section, $N(\rho,\tau )$ is one point. \ Hence the map
(3.4) becomes a map
\begin{equation*}
\phi :\frak{h}^{3}-\Delta \rightarrow G/T
\end{equation*}

compatible with $W\times SU(2).$ \ This is the result, generalizing the case
of $U(n),$ which arose from the Berry-Robbins paper and provided our original
motivation.

\item[(c)]  $\rho $ the sub-regular\footnote{ This means that the corresponding nilpotent orbit in ${\frak g}^C$ is subregular (the unique codimension $2$ orbit in the nilpotent variety).} $SU(2).$ \ Then, as we shall see later, $%
N(\rho ,\tau )$ is the $4$-dimensional ALE space studied by Kronheimer \cite{KronALE}. In this case $Z(\rho )$ is finite for all simple $G$ except $SU(n)$ when it is
$U(1).$ \ This circular symmetry corresponds to the Gibbons-Hawking
construction \cite{GH}.
\end{enumerate}

\section{The complex picture}

In this section we shall break the symmetry of ${\Bbb R}^{3}$ by picking a preferred
axis and consider the orthogonal projection
\begin{equation*}
\pi :{\Bbb R}^{3}\rightarrow {\Bbb R}^{2}\cong {\Bbb C},
\end{equation*}
identifying ${\Bbb R}^{2}$ with the complex plane. \ The symmetry group $SO(3)$ is
then reduced to $SO(2)=U(1).$

The preferred axis will pick out a distinguished complex symplectic structure
on all the hyperk\"{a}hler manifolds described in the preceding sections. \
We shall now analyse the complex manifolds that arise. \ In this we are
essentially following Kronheimer \cite{Kron2} as extended by the second author \cite{BielJLMS, BielAGAG}.

Let $\tau =(\tau _{1},\tau _{2},\tau _{3})$ be a regular triple with
projection $\pi (\tau )=\sigma =\tau _{2}+i\tau _{3}.$ \ We shall in the
first instance assume that $\sigma $ is a regular point of the complex Lie
algebra $\frak{h}\otimes {\Bbb C}$. Then the main result proved in \cite{BielJLMS}
identifies the preferred complex symplectic structure of the manifold $N(\rho
,\tau ).$ \ To describe this we need to recall the slice $S(\rho )$
introduced by Slodowy \cite{Slo}. \ First we extend $\rho $ to a homomorphism of complex Lie algebras
\begin{equation*}
\rho :{\frak sl}(2,C)\rightarrow \frak{g}\otimes C.
\end{equation*}

Let
\begin{equation*}
h=\left(
\begin{array}{rr}
-1 & 0 \\ 0 & 1
\end{array}
\right) \quad x=\left(
\begin{array}{rr}
0 & 1 \\ 0 & 0
\end{array}
\right) \quad y=\left(
\begin{array}{cc}
0 & 0 \\ 1 & 0
\end{array}
\right)
\end{equation*}
be the standard basis of $\frak{sl} (2,C)$ and let $H,X,Y$ be their images
under $\rho .$ \ We put
\begin{equation}
S(\rho )=Y+Z(X)  \label{4.1}
\end{equation}
where $Z(X)$ is the centralizer of $X$ in $\frak{g}\otimes C.$ \ Then $S(\rho
)$ is a transverse slice to the orbit of $Y.$ \ It is transverse
to any adjoint orbit of $G^{\Bbb C}$ it meets. \ In particular it is transverse to
the orbit $G^{\Bbb C}\sigma$ and so intersects this in a manifold. \ Then we
have \cite{BielJLMS}
\begin{equation}
N(\rho ,\tau )\cong G^{\Bbb C}\sigma \cap S(\rho )  \label{4.2}
\end{equation}
where $N(\rho ,\tau )$ is given its preferred complex structure. \ Varying $%
\tau _{1},$ while keeping $\sigma =\tau _{2}+i\tau _{3}$ fixed, gives
different K\"{a}hler metrics to the complex manifold in (4.2).

If $\tau _{1}$ is a regular point of $\frak{h},$ then $(\tau _{1},0,0)$ is a
regular triple so that $N(\rho ,\tau )$ is still a complex manifold for $%
\sigma =0,$ where the isomorphism (4.2) breaks down. \ To understand what
happens here we have to explain the Brieskorn-Grothendieck theory of the
simultaneous resolution.

Starting now with the complex Lie group $G^{\Bbb C}$ we let $B$ be a Borel
subgroup, $\frak{b}$ its Lie algebra and $\frak{h}^{\Bbb C}$ a Cartan subalgebra
in $\frak{b}.$ \ The Grothendieck resolution is then given by the diagram
\begin{equation}
\begin{array}{ccccc}
G^{\Bbb C}\times_{B}\frak{b} &  & \overset{\psi }{\rightarrow } &  & \frak{g}^{\Bbb C}
\\ &  & &  &
\\ \downarrow \theta &  &  &  & \downarrow \chi \\ &  &  &  &  \\ \frak{h}^{\Bbb C} &  &
\rightarrow &  & \frak{h}^{\Bbb C}/W
\end{array}
\label{4.3}
\end{equation}
where $B$ acts on $G^{\Bbb C}$ on the right and by the adjoint action on
$\frak{b}.$ \ The vertical maps are given by taking the semi-simple parts
(the
``eigenvalues''). \ The key property of this diagram is that the fibres of $%
\theta $ provide resolutions of the singularities of the fibres of $\chi $
and that $\theta $ is a smooth fibration (and topologically a product). \
Note in particular that $\theta ^{-1}(0)$ is a resolution of the nilpotent
variety $\frak{N}:$ it is isomorphic to the cotangent bundle $T^{\ast
}(G^{\Bbb C}/B).$

We can now restrict this diagram to the slice $S(\rho )$, giving the diagram
\begin{equation}
\begin{array}{ccccc}
\psi ^{-1}(S(\rho )) &  & \rightarrow  &  & S(\rho ) \\ &  &  &  &  \\
\downarrow \theta (\rho ) &  &  &  & \downarrow \chi (\rho ) \\ &  &  &  &
\\ \frak{h}^{\Bbb C} &  & \rightarrow  &  & \frak{h}^{\Bbb C}/W
\end{array}
\label{4.4}
\end{equation}
Again the fibres of $\theta (\rho )$ resolve the singularities of the fibres
of $\chi (\rho )$ and $\theta (\rho )$ is a smooth fibration. \ In particular
the inverse image $\theta (\rho )^{-1}(0)$ resolves the singularities of
$\frak{N}\;\cap S(\rho ).$

The generic fibre of $\chi (\rho )$ is the manifold $G^{\Bbb C}(\sigma )\cap
S(\rho )$ of (4.2). \ As shown in \cite{BielCMP, BielCMP2} the manifold $\psi
^{-1}S(\rho )$ of (4.4) can be naturally identified with the submanifold $%
N_{\tau _{1}}(\rho )\subset N(\rho )$ (with fixed regular $\tau _{1}$). \ In
other words the complex manifolds $N(\rho ,\tau )$ are the fibres of $\theta
(\rho )$ and in particular
\begin{equation}
N(\rho ;\tau _{1},0,0)  \label{4.5}
\end{equation}
\textbf{is the Grothendieck resolution of }$\frak{N}\;\cap S(\rho ).$

Let us illustrate all this by examining the three special cases of $\rho :$

\begin{enumerate}
\item[(a)]  $\rho =0,\;S(\rho )=\frak{g}^{\Bbb C},\;N(0;\tau )=G^{\Bbb C}(\sigma )$ and $%
N(0;\tau _{1},0,0)$ is the resolution of $\frak{N}$ and diagram (4.4) is just
(4.3).

\item[(b)]  $\rho $ the regular $\frak{su}(2),\;S(\rho )$ is a translate of $\frak{h}^{\Bbb C}/W,$
the manifold in (4.2) is just a point and $\theta (\rho )$ is an isomorphism.

\item[(c)]  $\rho $ the sub-regular $\frak{su}(2),$ the manifolds in (4.2) have
complex dimension $2$ and the fibres of $\theta (\rho )$ are the ALE spaces
as discussed by Kronheimer \cite{KronALE}.
\end{enumerate}

Considering again the general case, we have a fibration $N(\rho )\rightarrow
h^{3}-\Delta $ with hyperk\"{a}hler manifolds $N(\rho ,\tau )$ as fibres. \
The group $SU(2)$ acts on this fibration. \ The subgroup $U(1)=SO(2)\subset
SO(3)$ fixing a direction of ${\Bbb R}^{3}$ has fixed points of the form $\tau
=(\tau _{1},0,0)\in h^{3}-\Delta $ and so its double-cover $S\subset SU(2)$
acts on the fibre $N(\rho ,\tau )$ over this point. \ As we have seen this
fibre in its complex structure fixed by $S$ is a complex manifold which
can be identified with the resolution of $\frak{N}\;\cap S(\rho )$. \ Thus $%
\frak{N}$ $(\rho ,\tau )$ has a holomorphic action of $U(1)$, in addition to
a commuting action of $Z(\rho ).$ 
As this holomorphic action of $U(1)$ must leave $Y$ fixed, it is the composition of the complex scalar action on $\frak{g}^{\Bbb C}$ and of the adjoint action by $\rho(U(1))$. This will be explained more fully in section \ref{analytic}.
\ The map
\begin{equation*}
\frak{N}(\rho ,\tau )\rightarrow S(\rho )
\end{equation*}
defines a distinguished compact complex subspace which is the inverse image
of the base point $Y\in S(\rho )$ (see (4.1)). \ From the Grothendieck
resolution (4.4) we see that this is just \textbf{the fixed point set of the
action of }$Ad(Y)$ on $G^{\Bbb C}/B=G/T.$ \ Equivalently, viewing $G^{\Bbb C}/B$ as the
space $\frak{B}$ of all Borel subgroups, it is the set of all Borel
subgroups whose Lie algebra contain $Y.$ \ We shall denote it by $\frak{B}%
_{Y}.$ \ When $\rho $ is the regular $SO(2),$ $Y$ is regular and $\frak{B}%
_{Y}$ is a point. When $\rho $ is sub-regular, \ $Y$ is sub-regular and $%
\frak{B}_{Y}$ is $1$-dimensional, consisting of rational curves intersecting
as in the Dynkin diagram \cite{Slo}. \ In general $\frak{B}_{Y}\subset \frak{N}%
(\rho ,\tau )$ is the ``compact core'' of the open manifold, and carries all
its topological information. \ More precisely, \textbf{the action of }$S$
extends to an action of $C^{\ast }$ \textbf{all of whose orbits have limits
(as }$z\rightarrow \infty )$ in $\frak{B}_{Y}.$ \ The observation essentially
goes back to Slodowy \cite{Slo} and will be recalled in detail in the next section.

\section{Relation with Kazhdan-Lusztig}

In a long series of papers (see \cite{KL1, KL}) Kazhdan and Lusztig made an
extensive study of representation of the Hecke algebras $H$ associated to
Weyl groups (both finite and affine). \ A comprehensive account of this theory is given in \cite{CG}. These algebras are defined over the
finite Laurent series
\begin{equation*}
A={\Bbb C}[q,q^{-1}]
\end{equation*}
and reduce to the group algebras of the Weyl group when $q=1.$

Kazhdan and Lusztig construct representations of $H$ on the equivariant $K$%
-groups of certain subspaces of the flag manifold of the Lie group $G.$ \ The
purpose of this section is to show how all the ingredients in the
Kazhdan-Lusztig construction arise naturally in our context. \ It is our hope
that this will shed light on the geometric significance of the Hecke
algebras. \ Essentially, by using the ``quaternionic'' aspect of Lie groups
which we have been emphasizing we are able to move outside the purely complex
theory of Lie groups where Kazhdan and Lusztig work. \ Since they use the
Grothendieck resolution (and ideas of Brieskorn and Slodowy) it is not
surprising that the hyperk\"{a}hler story described in previous sections
should be relevant.

Given $\rho :SU(2)\rightarrow G$ we recall that we have the fibration (3.6)
\begin{equation*}
\begin{array}{c}
N(\rho ) \\
\\
\downarrow \\
\\
\frak{h}^{3}-\Delta
\end{array}
\end{equation*}
whose fibres $N(\rho ,\tau )$ are hyperk\"{a}hler manifolds, and that the
group
\begin{equation*}
W\times SU(2)\times Z(\rho )
\end{equation*}
acts on the fibration (where $Z(\rho )$ centralizes the image of $\rho ).$ \
We now fix a direction in ${\Bbb R}^{3}$ reducing the $SU(2)$ symmetry to a circle $%
S.$ \ We identify the ring $A$ with the character ring of $S$ (over ${\Bbb C})$%
\begin{equation}
A=R(S)\otimes {\Bbb C}  \label{5.1}
\end{equation}
(since $S\subset SU(2)$ double-covers $SO(2)\subset SO(3)$, our $q$ is the
square-root of the one in \cite{KL}.

This means that any space $X$ on which $S$ acts will have an equivariant $K$%
-group
\begin{equation}
K_{S}(X)\otimes {\Bbb C}  \label{5.2}
\end{equation}
which is an $A$-module. \ If $X$ is not compact we shall use $K$-theory with
\textbf{compact supports} in (5.2)

Consider now a fixed point $\tau $ for the action of $S$ on
$\frak{h}^{3}-\Delta .$ \ If we choose our coordinates of ${\Bbb R}^{3}$ so that $S$
defines rotation in the $(x_{2},x_{3})$ plane then $\tau $ is fixed under $S$
if it is of the form $(\tau _{1},0,0).$ \ Note that the set of such points
can be identified with the regular points of $h$ and so the components are
permuted by the Weyl group. \ A choice of component is essentially the same
as a choice of Borel subgroup of $G^{\Bbb C}$ containing $T,$ or equivalently a
choice of complex structure on the flag manifold $G/T.$

The fibre $N(\rho ,\tau )$ over $\tau $ has a complex structure (singled out
by our choice of direction) and a holomorphic action of $S.$ \ We can
therefore consider the $K$-group (with compact support)
\begin{equation}
K_{S}(N(\rho ,\tau ))\otimes C  \label{5.3}
\end{equation}

Inside $N(\rho ,\tau )$ we have its ``compact core'', namely the fixed-point
set $\frak{B}_{Y}$ of the nilpotent element $Y\in {\frak sl}(2, {\Bbb C}),$ as explained in
\S 4, and the action of ${\Bbb C}^{\ast }$ (complexification of $S$) has all limits
$z\rightarrow \infty $ in $\frak{B}_{Y}.$

Now Kazhdan-Lusztig work with the ``homology'' version $K_{0}^{S}$ of $%
K_{S}^{0}$ and observe that, in the situation just described, we have a
natural isomorphism
\begin{equation}
K_{0}^{S}(\frak{B}_{Y})\otimes C\cong K_{S}^{0}(N(\rho ,\tau ))\otimes {\Bbb C}.
\label{5.4}
\end{equation}

It is modules such as these in (5.4) (and various refinements) that are the $%
A$-modules studied by Kazhdan and Lusztig. \ One obvious refinement is to
enhance the symmetry from $S$ to $S\times Z(\rho ),$ or to a subgroup of
this.

The Weyl group $W$ does not act on the spaces in (5.4), it permutes them. \
However we also have the map
\begin{equation*}
\phi :N(\rho )\rightarrow G/T
\end{equation*}
defined by (3.4) and this is compatible with the action of $W\times Z(\rho ).$ \
This makes the groups in (5.4) into modules over
\begin{equation*}
K_{S}(G/T)\otimes {\Bbb C}
\end{equation*}
(where $S$ acts on $G/T$ via $\rho )$ and more generally we can replace $S$
by $S\times Z(\rho ).$

Let us now describe why this picture might help to explain the geometric
significance of the Hecke algebra and its modules. \ As we have seen the $K$%
-groups in question, disregarding for the moment the $S$-equivariance, are $K
$-groups of fibres over $\frak{h}^{3}-\Delta $ with an action of $W$ on the fibration. \
Alternatively they are $K$-groups of fibres over $(\frak{h}^{3}-\Delta )/W.$
\ In a non-equivariant situation this gives rise to the monodromy action of
$W.$ \ The action of $W$ on the homology of the fibres essentially gives the
Springer representations. \ In an equivariant situation (e.g. with an $S$%
-action) it is not clear what replaces monodromy, since $S$ only acts on
fibres over its fixed points. \ This suggests that $(\frak{h}^{3}-\Delta
)/W,$ together with its $S$-action, somehow produces the Hecke algebra
(instead of the fundamental group) and that bundles over this space (together
with compatible $S$-action) yield $H$-modules. \ One small piece of evidence
in favour of this idea is to note that the $S$-equivariant analogue of a path
from a point $\tau $ (fixed by $S)$ to its transform $\omega (\tau
),\;\omega \in W,$ is a $2$-sphere acted on by $S.$ \ The equivariant $K$%
-theory of such a $2$-sphere is an $A$-module with one generator, satisfying
a quadratic equation which is essentially the defining equation for generators
of $H.$

Unfortunately, although this is an appealing idea, we have not yet seen how
to carry it out. \ What we have done however is to put the general Kazhdan-Lusztig
construction into a more natural form.

In particular the circle symmetry is enlarged to a full $SU(2)$-action. \ We
hope the pay-off will emerge later.

\section{Analytic details\label{analytic}}

In this section we shall explain the analytic details behind the main
construction and in particular show how to define the spaces $N(\rho)$ and
$M(\rho)$ of section 3 as moduli spaces of solutions to Nahm's equations. The
Nahm equations will be the full translation-invariant anti-self-duality
equations on ${\Bbb R}^4$:
\begin{equation}\dot{T}_i+[T_0,T_i]+[T_k,T_j]=0\;,\label{Nahm}\end{equation}
where $(i,j,k)$ run over cyclic permuations of $(1,2,3)$. This form of Nahm's equations admits an action by the gauge group of $G$-valued functions $g(t)$:
\begin{eqnarray} T_0&\mapsto & \Ad(g)T_0-\dot{g}g^{-1}\nonumber\\ T_i&\mapsto & \Ad(g)T_i\;,\;\;\qquad i=1,2,3.\label{action}\end{eqnarray}
The component $T_0$ can be gauged away if we allow arbitrary gauge transformations.
We recall that the space $N^\prime(\rho)$ was defined as the space of
solutions to Nahm's equations on the half-line with poles of type $\rho$ at
$t=0$ and approaching a regular commuting triple as $t\rightarrow +\infty$.
As Kronheimer \cite{Kron} observes such a solution must approach its limit
exponentially fast.
\par
Let  $\Omega$ be the space of exponentially fast decaying functions in
$C^1[0,+\infty]$, i.e.:
\begin{equation}\Omega=\left\{f:(0,\infty]\rightarrow {\frak g}; \exists_{\eta >0}\sup_{t\geq 0}\left(e^{\eta t}\|f(t)\|+e^{\eta t}\|df/dt\|\right)<+\infty\right\}.\label{omega}\end{equation}
To define $N(\rho)$ let us fix a Cartan subalgebra $\frak{h}$ of $\frak{g}$
and consider solutions to \eqref{Nahm} on the half-line satisfying the
following boundary conditions at infinity:
\begin{itemize}
\item[(i)] $T_0(+\infty)=0$;
\item[(ii)]  $T_i(+\infty)\in {\frak h}$  for $i=0,\ldots,3$;
\item[(iii)] $(T_1(+\infty),T_2(+\infty),T_3(+\infty))$ is a regular triple, i.e. its centralizer is ${\frak h}$;
\item[(iv)]  $\left(T_i(t)-T_i(+\infty)\right)\in \Omega$ for $i=0,1,2,3$.\end{itemize}

In addition, the boundary conditions at $t=0$ are the same as for $N^\prime(\rho)$.
This space is acted upon by the gauge group ${\cal G}$ whose Lie algebra
consists of  {\em bounded} $C^2$-paths $\rho:[0,+\infty)\rightarrow {\frak
g}$ with $\rho(0)=0$ and $\dot{\rho}, [\tau,\rho]$ both belonging to $
\Omega$ for any regular element $\tau$ of ${\frak h}$. This means that any
element of ${\cal G}$ is asymptotic to an element of $T=\exp\frak{h}$.
Observe that we have a free action of $W=N(T)/T$ on $N(\rho)$ given by gauge
transformations asymptotic to elements of $N(T)$.
\par
We claim that the moduli space we obtain is the space $N(\rho)$ defined by
the diagram \eqref{3.3}. Indeed, we see that we can always make $T_0$
identically zero via a gauge transformation $g(t)$ with $g(0)=1$. This gives
us a projection $N(\rho)\rightarrow N^\prime(\rho)$. Now suppose we have two
solutions $(T_i)$ and $(T^\prime_i)$ in $N(\rho)$ which map to the same
element of $N^\prime(\rho)$. This means that $(T_i)$ and $(T^\prime_i)$ are
gauge equivalent via a gauge transformation $g(t)$ with $g(0)=1$. Moreover,
as the limit of both $(T_i)$ and $(T^\prime_i)$ is a regular triple in the
same Cartan subalgebra, $g(t)$ is asymptotic to an element of $N(T)$ and so
$(T_i)$ and $(T^\prime_i)$ are in the same $W$-orbit.
\par
The manifold $N(\rho)$ is not a hyperk\"ahler. Nevertheless it is fibred by the hyperk\"ahler manifolds $N(\rho,\tau)$ defined by fixing the limit $\tau=(\tau_1,\tau_2,\tau_3)$ of $T_1,T_2,T_3$ \cite{BielJLMS} (this is the fibration defined in \eqref{3.4}). As pointed out in section 2, a moduli space of solutions to Nahm's equations is expected to carry a hyperk\"ahler structure if it can be (formally) realised as a an infinite-dimensional hyperk\"ahler quotient. The spaces $N(\rho,\tau)$ are such quotients of the flat affine manifold consisting of all functions $(T_0,T_1,T_2,T_3)$ with prescribed boundary conditions.
\par
When $\sigma=\tau_2+i\tau_3$ is a regular element of ${\frak g}^{\Bbb C}$,  $N(\rho,\tau)$ has the complex structure (corresponding to choosing the $x_1$-axis in ${\Bbb R}^3$) described in \eqref{4.2}. In general, a complex structure of a hyperk\"ahler moduli space of solutions to Nahm's equations 
 can be identified  by writing Nahm's equations as equations for ${\frak g}^{\Bbb C}$-valued functions. If we choose an isomorphism (compatible with the usual metrics) $\,{\Bbb R}^3={\Bbb R}\times{\Bbb C}\,$, i.e. we choose  complex coordinates, say $\,(t+ix_1,x_2+ix_3)$, on ${\Bbb R}^4$, we can put
$$\alpha:=T_0+iT_1\,,\;\beta:=T_2+iT_3$$
The Nahm equations can then be written as:
\begin{equation}\frac{d\,}{dt}(\alpha+\alpha^{\ast})+[\alpha,\alpha^{\ast}]+[\beta,
\beta^{\ast}]=0\label{real}\end{equation}
and \begin{equation}\frac{d\,}{dt}\beta=[\beta,\alpha]\label{complex}\end{equation}
The second equation is preserved by the complex gauge transformations and our moduli space as a complex (in fact complex-symplectic) manifold is just $$\bigl(\text{solutions to \eqref{complex}}\bigr)/ \bigl( \text{complex gauge transformations}\bigr).$$ 
This is an example of identifying hyperk\"ahler and complex symplectic quotients \cite{HKLR}.   
\par
Returning to $N(\rho,\tau)$, we first observe, after Kronheimer \cite{Kron}, that when $\rho=0$, $N(\rho,\tau)$ is the complex adjoint orbit of $\sigma$ with the holomorphic identification given by 
\begin{equation} (\alpha(t),\beta(t))\mapsto \beta(0).\end{equation}
For a general $\rho$, $N(\rho,\tau)$ can be defined as the hyperk\"ahler quotient of the product manifold $N(0,\tau)\times N_\rho$, where $N_\rho$ is the moduli space of solutions to Nahm's equations on the interval $(0,1]$ with poles of type $\rho$ at $t=0$ and regular at $t=1$ (mod gauge transformations which are $1$ at both endpoints). This hyperk\"ahler manifold has been studied in detail in \cite{BielJLMS} where it was shown that with respect to any complex structure it is $S(\rho)\times G^{\Bbb C}$ ($S(\rho)$ is the transversal slice defined in \eqref{4.1}). In particular, when $\rho=0$, $N_\rho$ is isomorphic to $T^\ast G^{\Bbb C}$ as a complex-symplectic manifold \cite{Kron3,BielJLMS}.\newline
Both $N(0,\tau)$ and $N_\rho$ admit a hyperk\"ahler $G$-action given by gauge transformations with arbitrary values at $t=0$ and $t=1$, respectively. Taking the hyperk\"ahler quotient of $N(0,\tau)\times N_\rho$ by the diagonal $G$ is equivalent to gluing the solutions in $N_\rho$ at $t=1$ to those in $N(0,\tau)$ at $t=0$, and so it results in the manifold $N(\rho,\tau)$. On the other hand, the complex symplectic quotient of $\bigl(S(\rho)\times G^{\Bbb C}\bigr)\times {\cal O}(\sigma)$ by $G^{\Bbb C}$ is easily seen to be $S(\rho)\cap {\cal O}(\sigma)$ (the complex moment map on $S(\rho)\times G^{\Bbb C}$ is $\mu(\beta,g)=\Ad(g)\beta$ and on ${\cal O}(\sigma)$ it is the identity). The general mantra of identifying hyperk\"ahler and complex-symplectic quotients gives\footnote{ Strictly speaking a hyperk\"ahler quotient can in general only be identified with an open subset of a complex-symplectic quotient (of semi-stable points). The analytic argument that in our case the two coincide is given in \cite{BielJLMS}.} us the complex structure of $N(\rho,\tau)$.
\par
 $N(\rho)$ admits an action of $SU(2)$ defined as follows. Let $A$ be an element of $SU(2)$. Then ${A}$ acts on $N(\rho)$ by rotating the ``vector" $(T_1(t),T_2(t),T_3(t)$ and then acting on the resulting solution to Nahm's equations with a gauge transformation equal to $\rho({A})^{-1}$ at $t=0$.  This action leaves invariant the residues of $(T_0,T_1,T_2,T_3)$ at $t=0$.

\par
We shall now explain the diagram \eqref{3.5} in terms of the solutions to Nahm's equations. 
We shall define a torus bundle $M(\rho)$ over $N(\rho)$ which will be a
hyperk\"ahler manifold (more exactly: a hypercomplex manifold with a
compatible symmetric form which is generically non-degenerate). To define
this torus bundle  we first observe that  $N(\rho)$ can be also defined as
${\cal A}/{\cal G}$, where ${\cal A}$ is defined by omitting the condition
(i)  on the solutions to Nahm's equations in the definition of $N(\rho)$ and
the gauge group  is enlarged to ${\cal G}$ consisting of paths $g(t)$
asymptotic to $\exp(ht+\lambda h)$ for some $h\in {\frak h}$ and $\lambda\in
{\Bbb R}$. In other words  the Lie algebra of ${\cal G}$ consists of
$C^2$-paths $\rho:[0,+\infty)\rightarrow {\frak g}$ such that
\begin{itemize}
\item[(i)] $\rho(0)=0$ and $\dot{\rho}$ has a limit in ${\frak h}$ at $+\infty$;
\item[(ii)] $(\dot{\rho}-\dot{\rho}(+\infty))\in \Omega$, and $[\tau,\rho]\in\Omega$ for any regular element $\tau\in{\frak h}$;
\end{itemize}
The torus bundle $M(\rho)$ over $N(\rho)$ is defined as the quotient ${\cal
A}/{\cal G}_0$, where ${\cal G}_0$ is defined as  ${\cal G}$ with the added
condition:\newline
(iii) $\lim_{t\rightarrow +\infty} (\rho(t)-t\dot{\rho}(+\infty))=0$.
\par
In other words, elements $g(t)$ of  ${\cal G}_0$ are asymptotic to $\exp(ht)$
for some $h\in {\frak h}$. It is clear that   ${\cal G}/{\cal
G}_0=\exp({\frak h})$ and therefore $M(\rho)$ is a torus bundle over
$N(\rho)$.
\par
We observe that this $M(\rho)$ is the one defined by the diagram \eqref{3.5}.
Indeed, we can make $T_0$ identically zero by a gauge transformation
asymptotic to $g\exp{(ht)}$ where $g\in G$ and $h\in\frak{h}$. Since we
quotient by ${\cal G}_0$ we obtain a well defined element of $G$ fitting into
the diagram \eqref{3.5}  (observe that in the above description of $N(\rho)$, we obtain a gauge transformation asymptotic to $g\exp{(ht)}$ but defined only up to the action of $T$).
\medskip

The hyperk\"ahler structure of $M(\rho)$ is part of the general story discussed in section \ref{2} and its existence is proved in detail in \cite{BielCMP}. In particular the hyperk\"ahler moment map for the action of $T$ on $M(\rho)$ is given by $(T_1(+\infty),T_2(+\infty),T_3(+\infty))$, and so the hyperk\"ahler quotients are the fibers $N(\rho,\tau)$ of the map \eqref{3.6}. There is an action of $SU(2)$ defined on $M(\rho)$ in exactly the same way as for $N(\rho)$.
This action rotates the complex structures of $M(\rho)$ which are therefore all equivalent. $M(\rho)$ as a complex manifold is discussed at length in \cite{BielCMP}. In the complex picture only the action of the $U(1)\subset SU(2)$ preserving the chosen complex structure is visible. This is the action described towards the end of \S 4.

\end{document}